\documentclass[12pt]{article}
\usepackage{amssymb}
\usepackage{amsmath}
\usepackage{amsthm}
\usepackage{color}
\usepackage{comment}
\usepackage{graphicx}

\begin{document}

\def\cP{\mathcal{P}}
\definecolor{myBlue}{rgb}{0,.4,.6}

\newcommand{\removableFootnote}[1]{}

\newtheorem{theorem}{Theorem}[section]
\newtheorem{lemma}[theorem]{Lemma}
\newtheorem{proposition}[theorem]{Proposition}





\title{The instantaneous local transition of a stable equilibrium to a chaotic attractor in piecewise-smooth systems of differential equations.}
\author{D.J.W.~Simpson\\\\
Institute of Fundamental Sciences\\
Massey University\\
Palmerston North\\
New Zealand}




\maketitle

\begin{abstract}


An attractor of a piecewise-smooth continuous system of differential equations
can bifurcate from a stable equilibrium to a more complicated invariant set
when it collides with a switching manifold under parameter variation.
Here numerical evidence is provided to show that this invariant set can be chaotic.
The transition occurs locally (in a neighbourhood of a point)
and instantaneously (for a single critical parameter value).
This phenomenon is illustrated for the normal form of a boundary equilibrium bifurcation in three dimensions
using parameter values adapted from of a piecewise-linear model of a chaotic electrical circuit.
The variation of a secondary parameter reveals a period-doubling cascade to chaos with windows of periodicity.
The dynamics is well approximated by a one-dimensional unimodal map which explains this bifurcation structure.
The robustness of the attractor is also investigated by studying the influence of nonlinear terms.

\end{abstract}

\section{Introduction}
\label{sec:intro}
\setcounter{equation}{0}



Piecewise-smooth (PWS) systems of differential equations are used in diverse areas to model phenomena
involving abrupt events or fast transitions. 
Boundaries of regions where a PWS system takes a particular smooth functional form are called switching manifolds.
On these manifolds the differential equations may be discontinuous.
In this case orbits may {\em slide} along the manifolds in a manner formulated by Filippov \cite{Fi88}.
This Letter concerns PWS systems that are instead continuous on switching manifolds.
Such systems provide useful models of
electrical circuits \cite{BaVe01,Ts03}, 
neuronal systems \cite{Co08,ToGe03}\removableFootnote{
Also: \cite{Mc70,LiSc00,ItMu94}.
},
and other biological systems \cite{AkBr05,JoKo99}, for example.

As the parameters of a PWS system are varied, an equilibrium can collide with a switching manifold.
This is referred to as a boundary equilibrium bifurcation (BEB).
BEBs are inherently local: one or more invariant sets grow out of a point.
However, in some instances the bifurcation induces a global change in the alignment of separatrices effecting a global bifurcation,
an example is the so-called canard super-explosion \cite{DeFr13}.
For PWS discontinuous systems, BEBs admit many possible geometric configurations
and until recently such bifurcations resisted classification for systems of only two dimensions \cite{GuSe11,Gl15c,HoHo16}.

\begin{figure}[t!]
\begin{center}
\setlength{\unitlength}{1cm}
\begin{picture}(8,6)
\put(0,0){\includegraphics[height=6cm]{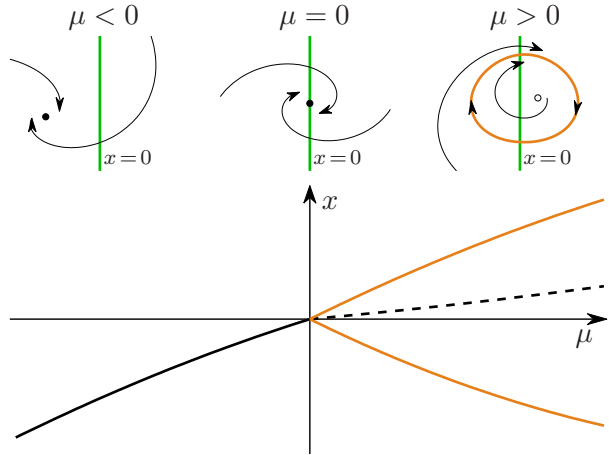}}
\put(7.52,1.51){\small $\mu$}
\put(4.15,3.27){\small $x$}
\put(.78,5.74){\small $\mu < 0$}
\put(3.58,5.74){\small $\mu = 0$}
\put(6.38,5.74){\small $\mu > 0$}
\put(1.25,3.87){\scriptsize $x \hspace{-.4mm} = \hspace{-.4mm} 0$}
\put(4.04,3.87){\scriptsize $x \hspace{-.4mm} = \hspace{-.4mm} 0$}
\put(6.83,3.87){\scriptsize $x \hspace{-.4mm} = \hspace{-.4mm} 0$}
\end{picture}
\caption{
A schematic bifurcation diagram of a supercritical Hopf-like BEB in a PWS continuous system.
As a parameter $\mu$ is increased through $0$, an equilibrium crosses a switching manifold ($x=0$).
Here the equilibrium loses stability (it is shown with a solid [dashed] line where it is stable [unstable])
and a stable periodic orbit is created (its minimum and maximum $x$-values are plotted).
Three representative phase portraits are also shown.
\label{fig:HopfLike}
}
\end{center}
\end{figure}

For PWS continuous systems, as an equilibrium crosses the switching manifold its associated eigenvalues change discontinuously
(and so it is sometimes called a discontinuous bifurcation).
One or more invariant sets may be created locally in the bifurcation.
All of these grow asymptotically linearly with respect to parameter change.
Some BEBs of PWS continuous systems mimic classical bifurcations, such as Hopf bifurcations \cite{FrPo97,SiMe07}, see Fig.~\ref{fig:HopfLike}. 
Others have no smooth analogue, such as the simultaneous creation of two equilibria and one periodic orbit \cite{LeNi04,Si10}.

For two-dimensional PWS continuous systems there are five types of non-degenerate BEBs
(this includes the case of an equilibrium crossing the switching manifold without
a change in stability or other invariant sets being created) \cite{FrPo98,SiMe12}.
No such classification is currently available in higher dimensions.
Recent studies of three-dimensional systems have revealed that the creation of periodic orbits
is typically governed by the nature of the dynamics on invariant cones \cite{CaFr05,CaFe12}.

BEBs can be analysed from a piecewise series expansion of the differential equations centred about the bifurcation.
Truncating the expansion to leading order generates a piecewise-linear (PWL) system with two pieces. 
This PWL approximation captures all structurally stable invariant sets created in the BEB.

Many authors have described chaotic dynamics in three-dimensional PWL continuous systems.
The most well-known example is probably Chua's circuit \cite{Ch94,Sh94},
but here the system has three pieces and so does not apply to BEBs.
The circuit system of Carroll \cite{Ca95,PiJa05}\removableFootnote{
Also \cite{Wu02}.
}
exhibits a chaotic Rossler-like attractor and is well modelled using only two pieces.
Shilnikov homoclinic chaos was described for a two-piece PWL continuous system by Sparrow in \cite{Sp81}.

The purpose of this Letter is to show that the local and instantaneous transition
from a stable equilibrium to a chaotic attractor
is possible in codimension-one BEBs of PWS continuous systems. 
In \S\ref{sec:3d} we describe the normal form for a BEB in three dimensions
in the case that a stable focus equilibrium exists on one side of the bifurcation and a saddle focus equilibrium exists on the other side.
In \S\ref{sec:returnMap} we construct a two-dimensional return map to explain the complicated oscillatory dynamics
on the side of the bifurcation with the saddle focus equilibrium.
In \S\ref{sec:RosslerLike} we show that the oscillatory dynamics can resemble a Rossler-like attractor
and identify a trapping region in the map.
Within this region the dynamics is well approximated by a one-dimensional unimodal map.
This explains the nature of bifurcation diagrams presented in \S\ref{sec:bifDiags}.
Here we also add nonlinear terms to the system
to illustrate how this type of BEB can be expected to behave in typical PWS systems.
Finally \S\ref{sec:conc} provides concluding remarks.

\section{BEBs in three dimensions with two foci}
\label{sec:3d}
\setcounter{equation}{0}

The normal form for a BEB of a PWS continuous system in three dimensions can be written as
\begin{equation}
\dot{X} = \begin{cases}
C_L X + e_3 \mu \;, & x \le 0 \\
C_R X + e_3 \mu \;, & x \ge 0
\end{cases} \;.
\label{eq:ocf}
\end{equation}
Here $X = (x,y,z) \in \mathbb{R}^3$ is the state variable
and $\mu \in \mathbb{R}$ is the parameter that governs the BEB which occurs at $\mu = 0$.
The switching manifold is the coordinate plane $x = 0$.
The matrices $C_L$ and $C_R$ are companion matrices
\begin{equation}
C_L = \begin{bmatrix}
\tau_L & 1 & 0 \\
-\sigma_L & 0 & 1 \\
\delta_L & 0 & 0
\end{bmatrix} \;, \qquad
C_R = \begin{bmatrix}
\tau_R & 1 & 0 \\
-\sigma_R & 0 & 1 \\
\delta_R & 0 & 0
\end{bmatrix} \;,
\label{eq:CLCR}
\end{equation}
where $\tau_L$, $\sigma_L$, and $\delta_L$
are the trace, second trace, and determinant of $C_L$,
and $\tau_R$, $\sigma_R$, and $\delta_R$ are the analogous quantities of $C_R$.
Throughout this Letter $e_1$, $e_2$ and $e_3$ denote the standard basis vectors of $\mathbb{R}^3$.
The form \eqref{eq:ocf}-\eqref{eq:CLCR} was first derived in a control theory context
and is also known as the observer canonical form \cite{CaFr02,DiMo11}\removableFootnote{
I might get away with omitting the paragraph:

For any constant $c > 0$, the system \eqref{eq:ocf} is invariant under the scaling $(X,\mu) \to (c X, c \mu)$.
Consequently, to describe the dynamics of \eqref{eq:ocf} it suffices to consider only $\mu \in \{-1,0,1\}$.
}.

If $\delta_L$ and $\delta_R$ are nonzero, \eqref{eq:ocf} has two potential equilibria
\begin{equation}
X_L = \frac{\mu}{\delta_L}
\begin{bmatrix} -1 \\ \tau_L \\ -\sigma_L \end{bmatrix} \;, \qquad
X_R = \frac{\mu}{\delta_R}
\begin{bmatrix} -1 \\ \tau_R \\ -\sigma_R \end{bmatrix} \;.
\nonumber
\end{equation}
The point $X_L$ is an equilibrium of \eqref{eq:ocf} if and only if $e_1^{\sf T} X_L \le 0$.
In this case we say $X_L$ is admissible; otherwise it is virtual.
Similarly $X_R$ is admissible if and only if $e_1^{\sf T} X_R \ge 0$.
At $\mu = 0$ both equilibria lie at the origin (here they are admissible and on the switching manifold).

The eigenvalues of $C_L$ and $C_R$ determine the stability of $X_L$ and $X_R$.
For the remainder of this Letter we assume these eigenvalues to be
\begin{equation}
\begin{split}
{\rm eig}(C_L):&~-\alpha_L \pm {\rm i} \beta_L,\;-\gamma_L, \\
{\rm eig}(C_R):&~\hspace{3.9mm}\alpha_R \pm {\rm i} \beta_R,\;-\gamma_R.
\end{split}
\label{eq:eigs}
\end{equation}
where $\alpha_L, \beta_L, \gamma_L, \alpha_R, \beta_R, \gamma_R > 0$.
We then have
\begin{equation}
\begin{aligned}
\tau_L &= -2 \alpha_L - \gamma_L \;, &
\tau_R &= 2 \alpha_R - \gamma_R \;, \\
\sigma_L &= \alpha_L^2 + \beta_L^2 + 2 \alpha_L \gamma_L \;, &
\sigma_R &= \alpha_R^2 + \beta_R^2 - 2 \alpha_R \gamma_R \;, \\
\delta_L &= -\left( \alpha_L^2 + \beta_L^2 \right) \gamma_L \;, &
\delta_R &= -\left( \alpha_R^2 + \beta_R^2 \right) \gamma_R \;,
\end{aligned}
\nonumber
\end{equation}
from which we see that
$X_L$ is an admissible stable focus for $\mu < 0$,
and $X_R$ is an admissible saddle focus for $\mu > 0$. 

\begin{figure}[b!]
\begin{center}
\setlength{\unitlength}{1cm}
\begin{picture}(8,6)
\put(0,0){\includegraphics[height=6cm]{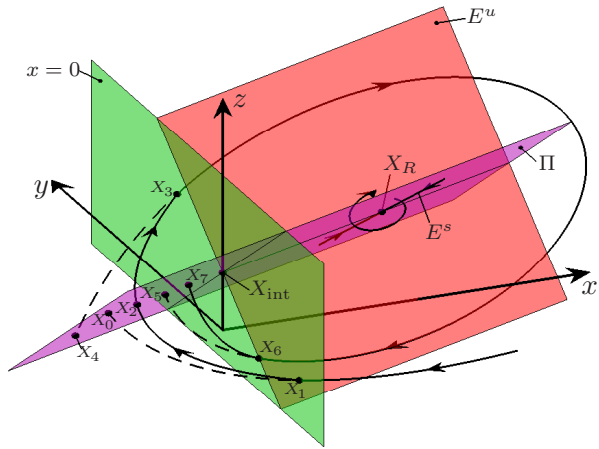}}
\put(7.76,2.1){\small $x$}
\put(.48,3.44){\small $y$}
\put(3.12,4.58){\small $z$}
\put(.38,5){\scriptsize $x=0$}
\put(6.23,5.7){\scriptsize $E^u$}
\put(5.68,2.83){\scriptsize $E^s$}
\put(7.19,3.75){\scriptsize $\Pi$}
\put(5.12,3.75){\scriptsize $X_R$}
\put(3.33,2.08){\scriptsize $X_{\rm int}$}
\put(1.23,1.68){\tiny $X_{\hspace{-.2mm}0}$}
\put(3.8,.76){\tiny $X_{\hspace{-.3mm}1}$}
\put(1.55,1.85){\tiny $X_{\hspace{-.25mm}2}$}
\put(2.02,3.45){\tiny $X_{\hspace{-.2mm}3}$}
\put(1.04,1.25){\tiny $X_{\hspace{-.2mm}4}$}
\put(1.88,2.04){\tiny $X_{\hspace{-.2mm}5}$}
\put(3.47,1.325){\tiny $X_{\hspace{-.2mm}6}$}
\put(2.49,2.29){\tiny $X_{\hspace{-.2mm}7}$}
\end{picture}
\caption{
A sketch of the phase space of \eqref{eq:ocf}-\eqref{eq:CLCR} with \eqref{eq:eigs}
showing intersections (both real and virtual)
of an orbit with the Poincar\'{e} section $\Pi$ and the switching manifold $x=0$.
Dashed lines show virtual extensions of the orbit defined by evolving
the right half-system of \eqref{eq:ocf} in the left half-space.
\label{fig:ppSchem}
}
\end{center}
\end{figure}

With $\mu > 0$, let us consider the stable and unstable manifolds of $X_R$.
Since \eqref{eq:ocf} is PWL, in a neighbourhood of $X_R$
the parts of these manifolds that emanate from $X_R$
coincide with the stable and unstable subspaces of $X_R$.
These subspaces, which we denote by $E^s$ and $E^u$ and are illustrated in Fig.~\ref{fig:ppSchem},
have directions given by the eigenvectors of $C_R$.
Specifically, the stable subspace $E^s$ is a line with direction
\begin{equation}
v = \begin{bmatrix}
1 \\ -2 \alpha_R \\ \alpha_R^2 + \beta_R^2
\end{bmatrix} \;,
\label{eq:v}
\end{equation}
and the unstable subspace $E^u$ is a plane normal to
\begin{equation}
w = \begin{bmatrix}
\gamma_R^2 \\ -\gamma_R \\ 1
\end{bmatrix} \;,
\label{eq:w}
\end{equation}
where $w^{\sf T}$ and $v$ are left and right eigenvectors of $C_R$
corresponding to the eigenvalue $-\gamma_R$\removableFootnote{
$E^u$ intersects $x=0$ on the line
\begin{equation}
z = \frac{1}{\gamma_R} + \gamma_R y \;.
\nonumber
\end{equation}
}.
Further away from $X_R$ the stable and unstable manifolds of $X_R$ curve
due to crossings with the switching manifold $x = 0$.

\section{A two-dimensional return map}
\label{sec:returnMap}
\setcounter{equation}{0}


Orbits of \eqref{eq:ocf} {\em graze} the switching manifold
(that is, intersect it tangentially) along the $z$-axis. 
The $z$-axis intersects the unstable subspace $E^u$ at
\begin{equation}
X_{\rm int} = \frac{\mu}{\gamma_R} e_3 \;.
\nonumber
\end{equation}
Therefore if an orbit on $E^u$ grazes the switching manifold, it does so at $X_{\rm int}$.

We define the Poincar\'{e} section\removableFootnote{
As a function of $x$ and $y$, $\Pi$ is given by 
\begin{align}
z &= \frac{1}{\gamma_R} \Big( \left( \left( \alpha_R^2 + \beta_R^2 \right)
\left( \gamma_R - 2 \alpha_R \right) - 4 \alpha_R^2 \gamma_R \right) x \nonumber \\
&\quad-
\left( \alpha_R^2 + \beta_R^2 + 2 \alpha_R \gamma_R \right) y + \mu \Big) \;.
\label{eq:zOnPi}
\end{align}
}
\begin{equation}
\Pi = \left\{ X_R + c_1 \big( X_R - X_{\rm int} \big) + c_2 v ~\middle|~ c_1, c_2 \in \mathbb{R} \right\} \;.
\label{eq:Pi}
\end{equation}
This is a particularly useful cross-section of phase space because it contains
the stable subspace $E^s$ and the important point $X_{\rm int}$\removableFootnote{
Since I don't end up with an explicit expression for the return map,
it is not easy to argue that a particular Poincar\'{e} section is best.
For any section the global map can be written by using $X(t) = {\rm e}^{t C_R} X_0$.
But with $E^s \subset \Pi$ the return time is simply $\frac{2 \pi}{\beta_R}$
leading to a simple form the global map.
}.

A typical orbit that we would like to describe using a return map on $\Pi$ is shown in Fig.~\ref{fig:ppSchem}.
This orbit evolves in the right half-space ($x > 0$) until it reaches the point $X_1$.
It then evolves in the left half-space ($x < 0$), intersecting $\Pi$ at $X_2$, until it reaches $X_3$.
Then it revolves around $X_R$ 
next intersecting $\Pi$ near $X_{\rm int}$ at $X_7$
(intersections with $\Pi$ far from $X_{\rm int}$, corresponding to large values of $x$, are ignored).

Rather than using the induced Poincar\'{e} map
(which provides $X_7$ as a function of $X_2$)
it is easier to work with a return map involving {\em virtual} intersection points by using a discontinuity map.
This concept was originally conceived by Nordmark \cite{No91}
and is an invaluable tool for understanding grazing bifurcations of periodic orbits in PWS systems \cite{DiBu01,FrNo00}.

The idea is to use the flow of the right half-space
to follow the orbit from its intersections with the switching manifold
into the left half-space until it intersects $\Pi$.
As indicated in Fig.~\ref{fig:ppSchem},
by evolving the orbit under the right half-flow forwards from $X_1$ and $X_6$
we generate intersections with $\Pi$ at $X_0$ and $X_5$,
and by evolving backwards from $X_3$ we generate an intersection at $X_4$.
We then define the return map $\cP : \Pi \to \Pi$ by $\cP(X_0) = X_5$.
Even though $X_0$ and $X_5$ are not points of the orbit (they are virtual intersection points),
$\cP$ captures the dynamics because it is conjugate to the Poincar\'{e} map.

We write
\begin{equation}
\cP = \cP_{\rm global} \circ \cP_{\rm disc} \;,
\nonumber
\end{equation}
where $X_5 = \cP_{\rm global}(X_4)$ describes one revolution about $X_R$ governed purely by the right half-flow,
and $X_4 = \cP_{\rm disc}(X_0)$ is the {\em discontinuity map} that corrects
for the difference between the left and right half-flows.
In cases for which the orbit intersects $\Pi$ without entering the left half-space,
$\cP_{\rm disc}$ is taken to be the identity map.

Since $E^s \subset \Pi$, the evolution time of $\cP_{\rm global}$
(i.e.~the time taken to go from $X_4$ to $X_5$)
is equal to $\frac{2 \pi}{\beta_R}$ for any $X_4 \in \Pi$.
Any $X_4 \in \Pi$ can be expressed uniquely in terms of the values of $c_1$ and $c_2$ in \eqref{eq:Pi}.
In terms of these $(c_1,c_2)$-coordinates, $\cP_{\rm global}$ is given simply by\removableFootnote{
By converting this to $(x,y)$-coordinates we obtain
the following description of the global map:
\begin{widetext}
\begin{equation}
\begin{aligned}
x &\mapsto \left( {\rm e}^{\frac{-2 \pi \gamma_R}{\beta_R}} +
\frac{2 \alpha_R}{\gamma_R} \left( {\rm e}^{\frac{2 \pi \alpha_R}{\beta_R}} -
{\rm e}^{\frac{-2 \pi \gamma_R}{\beta_R}} \right) \right) x +
\frac{1}{\gamma_R}  \left( {\rm e}^{\frac{2 \pi \alpha_R}{\beta_R}} -
{\rm e}^{\frac{-2 \pi \gamma_R}{\beta_R}} \right) y +
\frac{1 - {\rm e}^{\frac{2 \pi \alpha_R}{\gamma_R}}}
{\left( \alpha_R^2 + \beta_R^2 \right) \gamma_R} \;, \\
y &\mapsto \frac{2 \alpha_R}{\beta_R} \left( -2 \alpha_R + \gamma_R \right)
\left( {\rm e}^{\frac{2 \pi \alpha_R}{\beta_R}} -
{\rm e}^{\frac{-2 \pi \gamma_R}{\beta_R}} \right) x +
\left( {\rm e}^{\frac{2 \pi \alpha_R}{\beta_R}} -
\frac{2 \alpha_R}{\gamma_R} \left( {\rm e}^{\frac{2 \pi \alpha_R}{\beta_R}} -
{\rm e}^{\frac{-2 \pi \gamma_R}{\beta_R}} \right) \right) y +
\frac{\left( -2 \alpha_R + \gamma_R \right) \left( 1 - {\rm e}^{\frac{2 \pi \alpha_R}{\gamma_R}} \right)}
{\left( \alpha_R^2 + \beta_R^2 \right) \gamma_R} \;,
\end{aligned}
\end{equation}
\end{widetext}
}
\begin{equation}
c_1 \mapsto {\rm e}^{\frac{2 \pi \alpha_R}{\beta_R}} c_1 \;, \qquad
c_2 \mapsto {\rm e}^{\frac{-2 \pi \gamma_R}{\beta_R}} c_2 \;.
\label{eq:globalMapc1c2}
\end{equation}

The discontinuity map $\cP_{\rm disc}$ is comprised of three steps: $X_0 \to X_1 \to X_3 \to X_4$.
An explicit expression for the flow during each step is available because \eqref{eq:ocf} is PWL.
However, the evolution time of each step is given by the root of transcendental equation.
For this reason a complete explicit expression for $\cP_{\rm disc}$ is too complicated to provide here,
but it can be evaluated numerically without the use of an ODE solver
(this was done for the computations below).
Note that since the intersection of $\Pi$ with the switching manifold is not the grazing line (the $z$-axis),
it is possible for $X_2$ to occur before $X_1$ or after $X_3$.

\section{A Rossler-like attractor}
\label{sec:RosslerLike}
\setcounter{equation}{0}

\begin{figure}[b!]
\begin{center}
\setlength{\unitlength}{1cm}
\begin{picture}(8,6)
\put(0,0){\includegraphics[height=6cm]{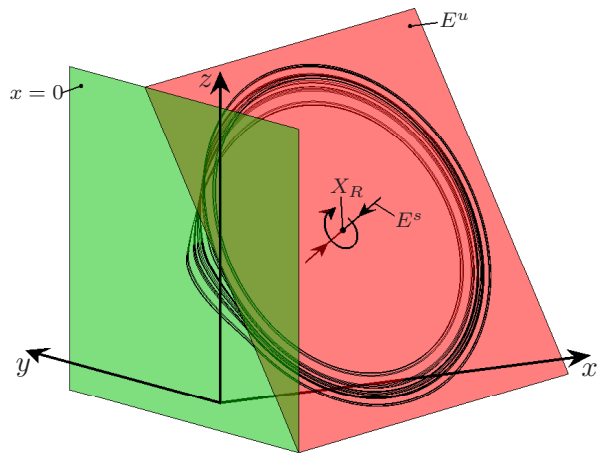}}
\put(7.6,1.08){\small $x$}
\put(.08,1.15){\small $y$}
\put(2.5,4.92){\small $z$}
\put(0,4.76){\scriptsize $x=0$}
\put(5.71,5.72){\scriptsize $E^u$}
\put(5.12,3.06){\scriptsize $E^s$}
\put(4.26,3.51){\scriptsize $X_R$}
\end{picture}
\caption{
The Rossler-like attractor of \eqref{eq:ocf} with \eqref{eq:exampleParams} and $\mu = 1$.
More precisely, part of an orbit (with transients removed) is shown.
The orbit was computed by iterating the return map $\cP$
and using the exact solution to each piece of \eqref{eq:ocf} to fill in the path of the orbit
between intersections with $\Pi$ and $x=0$.
\label{fig:RosslerLikeAttractor}
}
\end{center}
\end{figure}

Fig.~\ref{fig:RosslerLikeAttractor} shows the attractor of \eqref{eq:ocf}
with $\mu = 1$ and
\begin{equation}
\begin{aligned}
\alpha_L &= 0.3 \;, &
\beta_L &= 4 \;, &
\gamma_L &= 0.05 \;, \\
\alpha_R &= 0.02 \;, &
\beta_R &= 1 \;, &
\gamma_R &= 1 \;.
\end{aligned}
\label{eq:exampleParams}
\end{equation}
This example was obtained by putting Carroll's circuit model \cite{Ca95} in the normal form \eqref{eq:ocf}
and adjusting the parameters so that they are all positive and the attractor appears to remain chaotic.

Let us consider an orbit following the attractor.
While in the right half-space, the orbit approaches $E^u$ but spirals away from $X_R$.
While in the left half-space, the orbit is pulled away from $E^u$ but brought closer to $X_R$.
The orbit thus repeatedly lifts away from and folds back onto $E^u$
in a manner similar to the Rossler attractor \cite{Ro76}.

In view of the contraction towards $E^u$,
when the orbit intersects the switching manifold from the right half-space
it does so at points relatively close to $E^u$.
Therefore, for this orbit, points in the domain of $\cP$
(obtained by evolving the right half-flow until intersecting $\Pi$
even if this involves entering the left half-space)
are near $E^u$.
We write $\Pi \cap E^u = \{ g(x) \,|\, x \in \mathbb{R} \}$, where 
\begin{equation}
g(x) = \begin{bmatrix} x \\
(\gamma_R - 2 \alpha_R) x \\
-2 \alpha_R (2 \alpha_R + 1) x + \mu
\end{bmatrix} \;,
\label{eq:g}
\end{equation}
which follows from the above formulas for $\Pi$ and $E^u$.

\begin{figure}[t!]
\begin{center}
\setlength{\unitlength}{1cm}
\begin{picture}(8,6)
\put(0,0){\includegraphics[height=6cm]{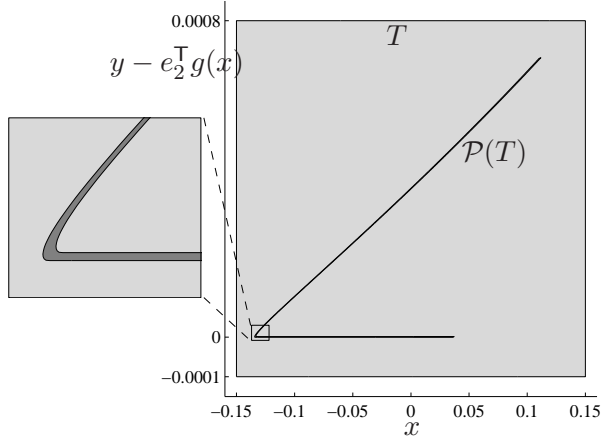}}
\put(5.33,.16){\small $x$}
\put(1.42,5){\small $y - e_2^{\sf T} g(x)$}
\put(5.1,5.35){\small $T$}
\put(6.1,3.8){\small $\cP(T)$}
\end{picture}
\caption{
A trapping region $T$ and its image under the return map $\cP$
for the parameter values \eqref{eq:exampleParams} and $\mu = 1$.
\label{fig:trappingRegion}
}
\end{center}
\end{figure}

Fig.~\ref{fig:trappingRegion} illustrates the action of $\cP$
on a rectangle $T$ for the same parameter values as Fig.~\ref{fig:RosslerLikeAttractor}.
The quantity $y - e_2^{\sf T} g(x)$ is plotted on the vertical axis.
This is a measure of the distance from $E^u$;
the horizontal axis thus represents $\Pi \cap E^u$.
Notice $\cP(T) \subset {\rm interior}(T)$,
thus $T$ is a trapping region and contains the attracting set $\bigcap_{i=0}^\infty \cP^i(T)$.
The vertical axis spans a relatively small range of values
which is a consequence of the contraction to $E^u$ described above.

\begin{figure}[t!]
\begin{center}
\setlength{\unitlength}{1cm}
\begin{picture}(8,6)
\put(0,0){\includegraphics[height=6cm]{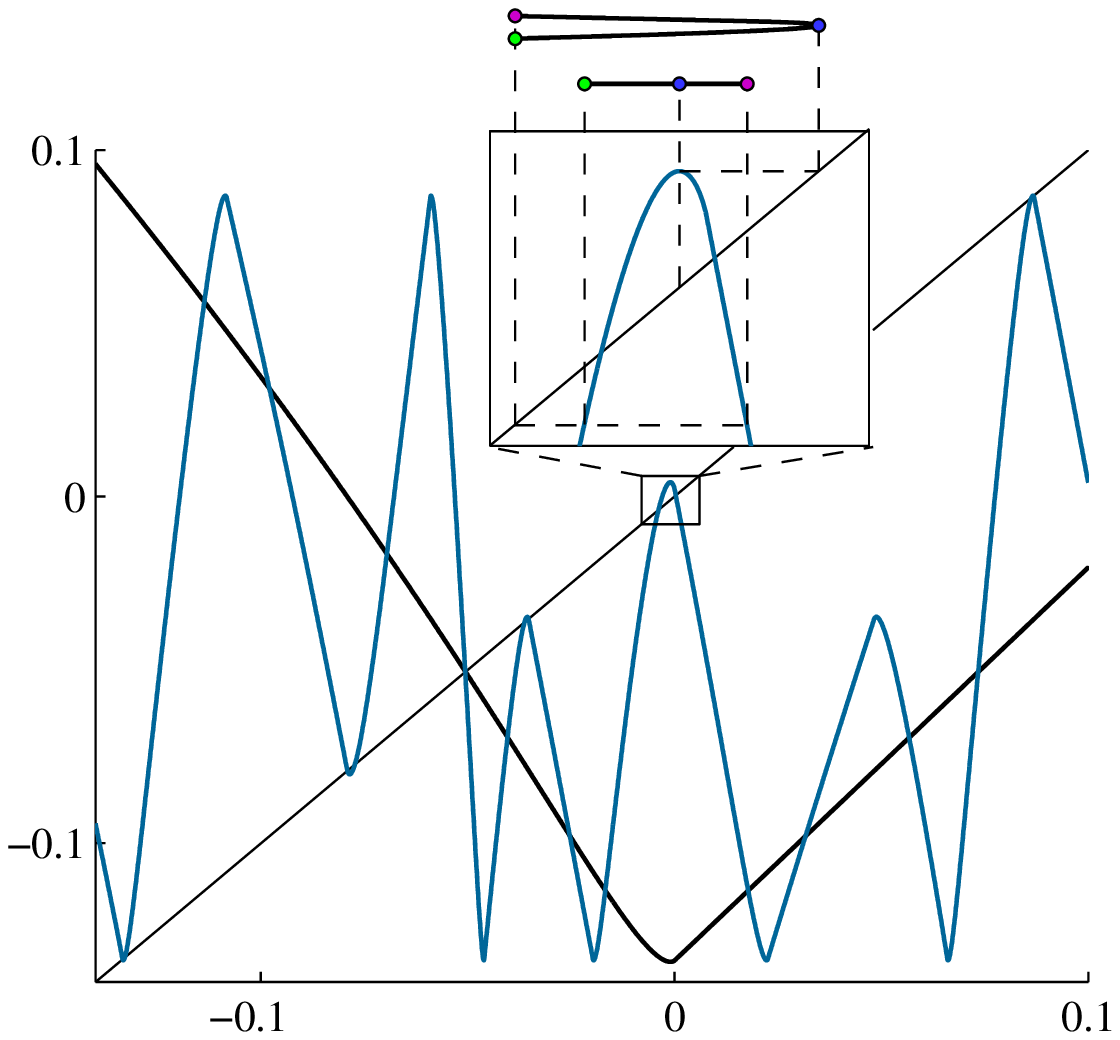}}
\put(4.34,.1){\small $x$}
\put(6.5,2.5){\footnotesize $f(x)$}
\put(6.58,4.3){\footnotesize \color{myBlue} $f^5(x)$}
\put(4.92,5.43){\scriptsize $I$}
\put(5.32,5.77){\scriptsize $f^5(I)$}
\end{picture}
\caption{
The one-dimensional map $f$ \eqref{eq:f} and its fifth iterate $f^5$
for the parameter values \eqref{eq:exampleParams} and $\mu = 1$.
\label{fig:oneDimMap}
}
\end{center}
\end{figure}

As a further benefit of the contraction to $E^u$, the one-dimensional map
\begin{equation}
f(x) = e_1^{\sf T} \cP(g(x)) \;,
\label{eq:f}
\end{equation}
captures the dynamics of $\cP$ with reasonable accuracy.
The map $f$ evaluates $\cP$ on $\Pi \cap E^u$ and retains only the $x$-component.
As shown in Fig.~\ref{fig:oneDimMap}, $f$ map is differentiable and unimodal.

The fifth iterate $f^5$ is also shown in Fig.~\ref{fig:oneDimMap}.
This has critical points near the $45^\circ$-line indicating that
the parameters \eqref{eq:exampleParams} are near that of a saddle-node bifurcation of a period-$5$ solution.
This observation helps us identify a Smale horseshoe
on which non-wandering dynamics is characterised by a shift on two symbols \cite{Ro04,Wi03}.
The image $f^5(I)$ of the interval $I = [-0.0045,0.0015]$
folds over $I$ and extends beyond both end points of $I$,
as illustrated in Fig.~\ref{fig:oneDimMap}.
Consequently $I$ can be fattened into a rectangle $R \subset \Pi$ for which $f^5(R)$ overlays $R$ as a horseshoe.
We do not show such a horseshoe as it is difficult to clearly discern the topology of $f^5(R)$
because it is both curved and strongly contracted.
A formal demonstration for the existence of a horseshoe,
as done in \cite{LlPo07} with $\sigma_L = \sigma_R = 1$, is beyond the scope of this Letter.
It should be achievable via direct but finicky and lengthy calculations using exact expressions for $\cP_{\rm disc}$.

\section{Bifurcation diagrams}
\label{sec:bifDiags}
\setcounter{equation}{0}

An increase in the value of $\gamma_L$ (which was fixed at $0.05$ in the previous section)
has a stabilising effect on the system.
This is because $-\gamma_L$ is an eigenvalue of $C_L$,
so $\gamma_L$ represents the strength of attraction of an eigen-direction of the left half-flow.

\begin{figure}[t!]
\begin{center}
\setlength{\unitlength}{1cm}
\begin{picture}(8,6)
\put(0,0){\includegraphics[height=6cm]{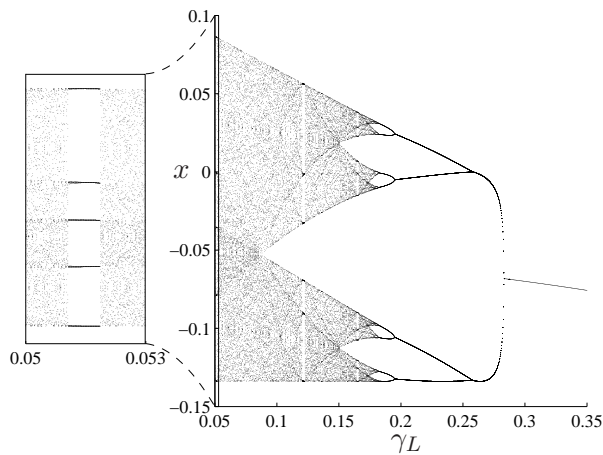}}
\put(5.13,.14){\small $\gamma_L$}
\put(2.24,3.72){\small $x$}
\end{picture}
\caption{
A bifurcation diagram for \eqref{eq:ocf} with \eqref{eq:exampleParams} and $\mu = 1$,
except that the value of $\gamma_L$ is variable.
The vertical axis shows $x$-values in the domain of $\cP$.
\label{fig:bifDiag_gammaL}
}
\end{center}
\end{figure}

Indeed with $\gamma_L = 0.35$ (and all other values unchanged from \eqref{eq:exampleParams})
there exists a stable periodic orbit.
Decreasing $\gamma_L$ from $0.35$ to $0.05$ reveals a period-doubling cascade to chaos, Fig.~\ref{fig:bifDiag_gammaL}.
The bifurcation diagram resembles that of the logistic map.
This is because the dynamics is well approximated by the one-dimensional unimodal map $f$.
Some windows of periodicity are visible;
a period-$5$ windows exists near $\gamma_L = 0.05$ as predicted above from the nature of $f^5$.

\begin{figure}[t!]
\begin{center}
\setlength{\unitlength}{1cm}
\begin{picture}(8,6)
\put(0,0){\includegraphics[height=6cm]{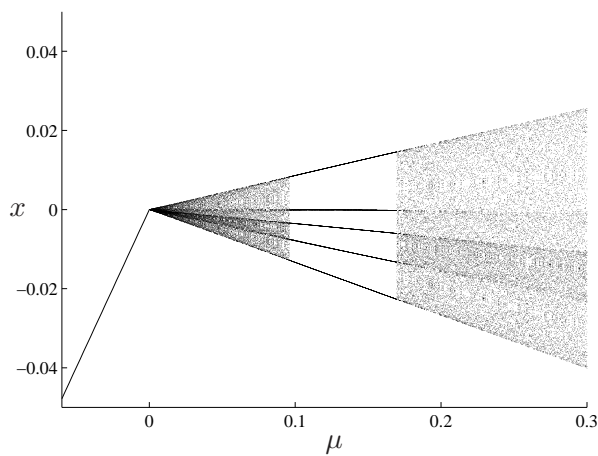}}
\put(4.3,.14){\small $\mu$}
\put(.1,3.225){\small $x$}
\end{picture}
\caption{
A bifurcation diagram for \eqref{eq:ocfNonlinear} with \eqref{eq:exampleParams}.
For $\mu < 0$, the $x$-value of the stable equilibrium is shown.
For $\mu > 0$, the $x$-values shown are those in the domain of $\cP$.
\label{fig:bifDiag_mu}
}
\end{center}
\end{figure}

In general, as parameters are varied to move away from a BEB,
the nonlinear terms that were omitted to produce the PWL approximation \eqref{eq:ocf}
have an increasingly greater influence.
To illustrate this influence on the creation of the Rossler-like attractor of Fig.~\ref{fig:RosslerLikeAttractor},
we simulate the system in the presence of an additional quadratic term:
\begin{equation}
\dot{X} = \begin{cases}
C_L X + e_3 \mu + e_1 x y \;, & x \le 0 \\
C_R X + e_3 \mu \;, & x \ge 0
\end{cases} \;.
\label{eq:ocfNonlinear}
\end{equation}
Fig.~\ref{fig:bifDiag_mu} shows the result of varying the value of $\mu$.
A period-$5$ window is visible; our numerical computations detected no stable periodic solutions for $\mu > 0$ to the left of this window.
It remains to determine if the periodic windows in
Figs.~\ref{fig:bifDiag_gammaL} and \ref{fig:bifDiag_mu} occur densely,
like for the logistic map \cite{GrSw97,Ly97},
or if chaos is robust over intervals of parameter values
which has been demonstrated in unimodal maps \cite{AnAl01}.

\section{Discussion}
\label{sec:conc}
\setcounter{equation}{0}

This paper shows how stable equilibria can bifurcate to chaotic attractors
in BEBs of PWS continuous systems.
The BEBs are local, not the result of suddenly accessing a distant part of phase space.
The analogous result for discrete-time systems has been known for some time \cite{DiGa98,YuBa98}.
More recently Glendinning has shown that stable fixed points
of PWS continuous maps can bifurcate to chaotic attractors of any positive integer dimension \cite{Gl15b}.

The dynamics described here is well-approximated by a one-dimensional unimodal map
and the existence of a Smale horseshoe was illustrated with numerical simulations.
A proof for the existence of a horseshoe seems attainable
via direct calculations of the return map $\cP$ because the differential equations are PWL.
In contrast, chaos in Chua's circuit was demonstrated using Shilnikov's theorem \cite{ChKo86}.

The phenomenon described here is perhaps best viewed as a one-parameter family of BEBs.
For each value of $\gamma_L$, the attractor created in the BEB at $\mu = 0$ is indicated by Fig.~\ref{fig:bifDiag_gammaL}.
As the value of $\mu$ is increased from $0$ there are additional bifurcations, as shown in Fig.~\ref{fig:bifDiag_mu} for $\gamma_L = 0.05$.
Both bifurcation diagrams show chaos with windows of periodicity.
It remains to determine if these windows occur densely
as this would establish whether or not the chaotic attractors are robust.


\begin{thebibliography}{10}

\bibitem{Fi88}
A.F. Filippov.
\newblock {\em Differential Equations with Discontinuous Righthand Sides.}
\newblock Kluwer Academic Publishers., Norwell, 1988.

\bibitem{BaVe01}
S.~Banerjee and G.C. Verghese, editors.
\newblock {\em Nonlinear Phenomena in Power Electronics.}
\newblock IEEE Press, New York, 2001.

\bibitem{Ts03}
C.K. Tse.
\newblock {\em Complex Behavior of Switching Power Converters.}
\newblock CRC Press, Boca Raton, FL, 2003.

\bibitem{Co08}
S.~Coombes.
\newblock Neuronal networks with gap junctions: {A} study of piecewise linear
  planar neuron models.
\newblock {\em SIAM J. Appl. Dyn. Sys.}, 7(3):1101--1129, 2008.

\bibitem{ToGe03}
A.~Tonnelier and W.~Gerstner.
\newblock Piecewise linear differential equations and integrate-and-fire
  neurons: Insights from two-dimensional membrane models.
\newblock {\em Phys. Rev. E}, 67:021908, 2003.

\bibitem{AkBr05}
O.E. Akman, D.S. Broomhead, R.V. Abadi, and R.A. Clement.
\newblock Eye movement instabilities and nystagmus can be predicted by a
  nonlinear dynamics model of the saccadic system.
\newblock {\em J. Math. Biol.}, 51:661--694, 2005.

\bibitem{JoKo99}
K.D. Jones and D.S. Kompala.
\newblock Cybernetic model of the growth dynamics of ${S}$\hspace{-.3mm}{\em
  accharomyces cerevisiae} in batch and continuous cultures.
\newblock {\em J. Biotech.}, 71:105--131, 1999.

\bibitem{DeFr13}
M.~Desroches, E.~Freire, S.J. Hogan, E.~Ponce, and P.~Thota.
\newblock Canards in piecewise-linear systems: explosions and superexplosions.
\newblock {\em Proc. R. Soc. A}, 469:20120603, 2013.

\bibitem{GuSe11}
M.~Guardia, T.M. Seara, and M.A. Teixeira.
\newblock Generic bifurcations of low codimension of planar {F}ilippov systems.
\newblock {\em J. Differential Equations}, 250:1967--2023, 2011.

\bibitem{Gl15c}
P.~Glendinning.
\newblock Classification of boundary equilibrium bifurcations in planar
  {F}ilippov systems.
\newblock {\em Unpublished.}, 2015.

\bibitem{HoHo16}
S.J. Hogan, M.E. Homer, M.R. Jeffrey, and R.~Szalai.
\newblock Piecewise smooth dynamical systems theory: the case of the missing
  boundary equilibrium bifurcations.
\newblock {\em Unpublished.}, 2016.

\bibitem{FrPo97}
E.~Freire, E.~Ponce, and F.~Torres.
\newblock {H}opf-like bifurcations in planar piecewise linear systems.
\newblock {\em Publicacions Matem\'{a}tiques}, 41:131--148, 1997.

\bibitem{SiMe07}
D.J.W. Simpson and J.D. Meiss.
\newblock Andronov-{H}opf bifurcations in planar, piecewise-smooth, continuous
  flows.
\newblock {\em Phys. Lett. A}, 371(3):213--220, 2007.

\bibitem{LeNi04}
R.I. Leine and H.~Nijmeijer.
\newblock {\em Dynamics and Bifurcations of Non-smooth Mechanical Systems},
  volume~18 of {\em Lecture Notes in Applied and Computational Mathematics}.
\newblock Springer-Verlag, Berlin, 2004.

\bibitem{Si10}
D.J.W. Simpson.
\newblock {\em Bifurcations in Piecewise-Smooth Continuous Systems.}, volume~70
  of {\em Nonlinear Science}.
\newblock World Scientific, Singapore, 2010.

\bibitem{FrPo98}
E.~Freire, E.~Ponce, F.~Rodrigo, and F.~Torres.
\newblock Bifurcation sets of continuous piecewise linear systems with two
  zones.
\newblock {\em Int. J. Bifurcation Chaos}, 8(11):2073--2097, 1998.

\bibitem{SiMe12}
D.J.W. Simpson and J.D. Meiss.
\newblock Aspects of bifurcation theory for piecewise-smooth, continuous
  systems.
\newblock {\em Phys. D}, 241(22):1861--1868, 2012.

\bibitem{CaFr05}
V.~Carmona, E.~Freire, E.~Ponce, and F.~Torres.
\newblock Bifurcation of invariant cones in piecewise linear homogeneous
  systems.
\newblock {\em Int. J. Bifurcation Chaos}, 15(8):2469--2484, 2005.

\bibitem{CaFe12}
V.~Carmona, S.~Fern\'{a}ndez-Garc\'{\i}a, and E.~Freire.
\newblock Saddle-node bifurcation of invariant cones in {3D} piecewise linear
  systems.
\newblock {\em Phys. D}, 241:623--635, 2012.

\bibitem{Ch94}
L.O. Chua.
\newblock Chua's circuit 10 years later.
\newblock {\em Internat. J. Circuit Theory Appl.}, 22(4):279--305, 1994.

\bibitem{Sh94}
L.P. Shil'nikov.
\newblock Chua's circuit: {R}igorous results and future problems.
\newblock {\em Int. J. Bifurcation Chaos}, 4(3):489--519, 1994.

\bibitem{Ca95}
T.L. Carroll.
\newblock A simple circuit for demonstrating regular and synchronized chaos.
\newblock {\em Am. J. Phys.}, 63:377--379, 1995.

\bibitem{PiJa05}
A.N. Pisarchik and R.~Jaimes-Re\'{a}tegui.
\newblock Homoclinic orbits in a piecewise linear {R}\"{o}ssler-like circuit.
\newblock {\em J. Phys.: Conf. Ser.}, 23:122--127, 2005.

\bibitem{Sp81}
C.~Sparrow.
\newblock Chaos in a three-dimensional single loop feedback system with a
  piecewise linear feedback function.
\newblock {\em J. Math. Anal. Appl.}, 83:275--291, 1981.

\bibitem{CaFr02}
V.~Carmona, E.~Freire, E.~Ponce, and F.~Torres.
\newblock On simplifying and classifying piecewise-linear systems.
\newblock {\em IEEE Trans. Circuits Systems I Fund. Theory Appl.},
  49(5):609--620, 2002.

\bibitem{DiMo11}
M.~di~Bernardo, U.~Montanaro, and S.~Santini.
\newblock Canonical forms of generic piecewise linear continuous systems.
\newblock {\em IEEE Trans. Automat. Contr.}, 56(8):1911--1915, 2011.

\bibitem{No91}
A.B. Nordmark.
\newblock Non-periodic motion caused by grazing incidence in impact
  oscillators.
\newblock {\em J. Sound Vib.}, 2:279--297, 1991.

\bibitem{DiBu01}
M.~di~Bernardo, C.J. Budd, and A.R. Champneys.
\newblock Normal form maps for grazing bifurcations in $n$-dimensional
  piecewise-smooth dynamical systems.
\newblock {\em Phys. D}, 160:222--254, 2001.

\bibitem{FrNo00}
M.H. Fredriksson and A.B. Nordmark.
\newblock On normal form calculation in impact oscillators.
\newblock {\em Proc. R. Soc. A}, 456:315--329, 2000.

\bibitem{Ro76}
O.E. R\:{o}ssler.
\newblock An equation for continuous chaos.
\newblock {\em Phys. Lett.}, 57A:397--398, 1976.

\bibitem{Ro04}
R.C. Robinson.
\newblock {\em An Introduction to Dynamical Systems. Continuous and Discrete.}
\newblock Prentice Hall, Upper Saddle River, NJ, 2004.

\bibitem{Wi03}
S.~Wiggins.
\newblock {\em Introduction to Applied Nonlinear Dynamical Systems and Chaos.},
  volume~2 of {\em Texts in Appl. Math.}
\newblock Springer-Verlag, New York, 2003.

\bibitem{LlPo07}
J.~Llibre, E.~Ponce, and A.E. Teruel.
\newblock Horseshoes near homoclinic orbits for piecewise linear differential
  systems in $\mathbb{R}^3$.
\newblock {\em Int. J. Bifurcation Chaos}, 17(4):1171--1184, 2007.

\bibitem{GrSw97}
J.~Graczyk and G.~Swiatek.
\newblock Generic hyperbolicity in the logistic family.
\newblock {\em Ann. Math.}, 146(1):1--52, 1997.

\bibitem{Ly97}
M.~Lyubich.
\newblock Dynamics of quadratic polynomials, {I-II}.
\newblock {\em Acta. Math.}, 178:185--297, 1997.

\bibitem{AnAl01}
M.~Andrecut and M.K. Ali.
\newblock Robust chaos in smooth unimodal maps.
\newblock {\em Phys. Rev. E}, 64(2):025203, 2001.

\bibitem{DiGa98}
M.~di~Bernardo, F.~Garofalo, L.~Glielmo, and F.~Vasca.
\newblock Switchings, bifurcations and chaos in {DC/DC} converters.
\newblock {\em IEEE Trans. Circuits Systems I Fund. Theory Appl.},
  45(2):133--141, 1998.

\bibitem{YuBa98}
G.~Yuan, S.~Banerjee, E.~Ott, and J.A. Yorke.
\newblock Border-collision bifurcations in the buck converter.
\newblock {\em IEEE Trans. Circuits Systems I Fund. Theory Appl.},
  45(7):707--716, 1998.

\bibitem{Gl15b}
P.~Glendinning.
\newblock Bifurcation from stable fixed point to ${N}$-dimensional attractor in
  the border collision normal form.
\newblock {\em Nonlinearity}, 28:3457--3464, 2015.

\bibitem{ChKo86}
L.O. Chua, M.~Komuro, and T.~Matsumoto.
\newblock The double scroll family.
\newblock {\em IEEE Trans. Circuits Systems}, CAS-33(11):1073--1118, 1986.

\end{thebibliography}
\end{document}